\documentclass[numbook, envcountsect, envcountsame, envcountreset, runningheads, smallextended]{svjour3}

\usepackage{amscd,amsfonts,amssymb,amsmath,latexsym,array,hhline,xcolor,graphicx}

\usepackage{latexsym}
\usepackage{times}

\newcommand\cE{{\cal E}}

\newcommand\cF{{\cal F}}

\newcommand\e{{\varepsilon}}

\def\E{{\bf E}}
\def\P{{\bf P}}

\def\text#1{\hbox{#1}}

\def\E{{\bf E}}
\def\P{{\bf P}}

\def\R{{\bf R}}

\def\build #1_#2{\mathrel{\mathop{\kern 0pt #1}\limits_{#2}}} 
\newcommand{\zs}[1]{{\mathchoice{#1}{#1}{\lower.25ex\hbox{$\scriptstyle#1$}}
{\lower0.25ex\hbox{$\scriptscriptstyle#1$}}}}

\numberwithin{equation}{section}

\newcommand\fdem{$\Box$}
\newcommand\beq{\begin{equation}}
\newcommand\eeq{\end{equation}}
\newcommand\bea{\begin{eqnarray}}
\newcommand\eea{\end{eqnarray}}
\newcommand\bean{\begin{eqnarray*}}
\newcommand\eean{\end{eqnarray*}}

\begin{document}

	\title{Ruin Probabilities for a Sparre Andersen Model with Investments:
the Case of Annuity Payments  
	}
	\author{ Yuri Kabanov \and Platon  Promyslov} 
	
	\institute{\at	Lomonosov Moscow State University, 
	Federal Research Center “Computer Science and Control” of the Russian Academy of Sciences
	Moscow, Russia, and  Universit\'e de Franche-Comt\'e, Laboratoire de Math\'ematiques, UMR CNRS 6623, 
	16 Route de Gray, 25030 Besan\c{c}on,  France \\
  \email{ykabanov@univ-fcomte.fr}. \\
	\and Lomonosov Moscow State University, Moscow, Russia\\
	 \email{platon.promyslov@gmail.com}}

%
%
%
%
%
\titlerunning{Ruin Probabilities for a Sparre Andersen Model with Investments}

\date{\today}
\maketitle
\qquad\qquad\qquad\qquad\qquad\qquad\qquad\qquad
{\sl Dedicated to the memory of Tomas Bj\"ork.}\\
\\

\begin{abstract}
This  note is a complement  to the paper by Eberlein, Kabanov, and Schmidt on the asymptotic of the ruin probability in a Sparre Andersen non-life insurance model  with investments a risky asset whose price follows a geometric L\'evy process.  
Using the techniques of semi-Markov processes   we extend the result of the mentioned paper to the case of annuities  and models with two-sided jumps.  
\end{abstract}

\smallskip

 


\section{Introduction}  
In the  classical Sparre Andersen model of  insurance company the counts of  claims form a renewal process.  
In  recent studies, see \cite{ACT}, \cite{EKT} and references therein,  this model was enriched  by the assumption that the capital reserve of the insurance company is  fully invested in a risky asset whose price evolves as a geometric L\'evy process. 
 In the paper \cite{EKT} by Eberlein, Kabanov, and Schmidt it was considered the non-life insurance version of such a model. 
It was shown that under rather mild 
hypotheses on the business process the asymptotic behavior of the (ultimate) ruin probability
 is essentially the same as in   the Cram\'er--Lundberg model with risky investments. Namely, the ruin probability decays, up to a multiplicative constant, as the function $u^{-\beta}$ where $u$, the 
 initial capital, tends to infinity. The decay rate $\beta$ depends only  of characteristics of the price process.  
 The method of analysis  in \cite{EKT} is based heavily on the assumption that the risk process has only downward jumps and, therefore, crosses the zero level only by a jump. This specific feature   allows a straightforward reduction to a discrete-time Markovian framework. 
 
 The approach of \cite{EKT} left an open question whether the results hold also in the case of   
upward jumps. This is a feature of the annuity model   when the risk process crosses the zero level in a continuous way. In a less popular mixed model with two-sided jumps the crossing may happen in both way.  
Of course, a positive answer is expected here:      
this was already established  for  the Cram\'er--Lundberg models with investments analyzed by Kabanov, Pergamenshchikov, and Pukhlyakov, \cite{KP}, \cite{KPukh},  as well as in  very general    L\'evy Ornstein--Uhlenbeck models introduced and studied by Paulsen, see \cite{Paul-93},  \cite{Paul},  \cite{Paul-98},  \cite{Paul-G}, and a more recent paper   \cite{KP2020} by Kabanov and Pergamenshchikov.  
 
Our note, based on the study \cite{EKT}, gives a positive answer  for a Sparre Andersen model with investments in its annuity version,  with the upward jumps.   We  discuss briefly  the needed changes leading to a result 
 for a model  with  upward and downward jumps used serving to describe the evolution of the capital reserve of a company with two types of the business activity.   
 Our  techniques is  based on the imbedding 
 of a semi-Markov process into  a Markov one by increasing the dimensionality. 
 
 In the paper we use standard notations of stochastic calculus and concepts discussed in details in 
 \cite{KP2020},  \cite{EKT}.

\section{The model} 
The Sparre Andersen model   with risky investments considered  contains two ingredients: 
\smallskip

1. The price process of a risky financial asset $S=(S_t)_{t\ge 0}$. It  is of  the form $S=\cE(R)$ where $\cE$ is the stochastic exponential,  $R$ is a L\'evy process with  the L\'evy triplet $(a,\sigma^2,\Pi)$  and such that  
$\Pi ((-\infty,-1])=0$. The latter condition ensures that the  jumps $\Delta R>-1$, hence,  the price $S>0$. In such a case, $S=e^V$ where  $V=\ln S$ is again a L\'evy process which can be given by the formula
\beq
\label{V}
V_t=at-\frac 12 \sigma^2 t + \sigma W_t+ h*(\mu-\nu)_t+(\ln (1+x)-h)*\mu_t, 
\eeq
where $h(x):=xI_{\{|x|\le 1\}}$. 
The L\'evy triplet of $V$  is  $(a_V,\sigma^2,\Pi_V)$ with 
$$
a_V=a-\frac{\sigma^2}{2}+\Pi(h(\ln(1+x))-h)  
$$ and 
$\Pi_V=\Pi\varphi^{-1}$,  $\varphi: x\mapsto \ln (1+x)$.  

It is assume that $R$ is non-deterministic, that is, at least one of the parameters $\sigma^2$ or $\Pi$
is not zero.

\smallskip
2. The "business process``. It is an independent of $S$ compound renewal process $P=(P_t)$.
 Classically, it can be written in  
the form 
\beq
\label{Pt}
P_t=ct+\sum_{i=1}^{N_t}\xi_i,  
\eeq 
where $N=(N_t)$ is a counting renewal process with the interarrival times (lengths of the inter jump intervals) $U_i:=T_i-T_{i-1}$, $i\ge 2$,  forming an i.i.d. sequence   independent of the i.i.d. sequence of random variables $\xi_i=\Delta P_{T_i}$, $i \ge 1$,  with the common law $F_\xi$, $F_\xi(\{0\})=0$. In the sequel, a ``generic`` r.v. with such a law is
denoted by $\xi$.  As usual, $T_0:=0$. The common law of  $U_i$ we denoted by $F$ and use the same character for its distribution function.  
\smallskip

%

The risk process  $X=X^{u}$, $u>0$,  is defined as the solution of the non-homogene\-ous linear stochastic  equation 
\beq
 \label{risk}
 X_t=u+  \int_0^t  X_{s-}  dR_s +P_t. 
 \eeq
 
 The ruin probability is the function of the initial capital  
 $\Psi(u):=\P[\tau^u<\infty]$ where $\tau^u:=\inf \{t: X^u_t\le 0\}$. 
 \smallskip

 The cases of major interest are: $c>0$ and $\xi_i<0$ (a non-life insurance model, considered in \cite{EKT}) and  
 $c<0$ and $\xi_i>0$ (annuities payments model). The latter case studied here, is often interpreted as a model of venture company paying salary and selling innovations. 
 The case where $F_\xi$ charges both half-axes 
 can be viewed as a model  of company combined two types of activity, see, e.g.,  \cite{ACT}. 
  and we study also the case where $\xi_i$ may take positive and negative values. 

If  $c\ge 0$ and $\xi>0$ the ruin never happens and this case is excluded from considerations.

 \medskip
 {\bf Standing assumption.} The cumulant generating function $H:q\to \ln  \E\,e^{-qV_{T_1}}$ of the random variable $V_{T_1}$ has a  root $\beta>0$ not laying on the boundary of the effective domain of $H$.  That is, if  the  ${\rm int\, dom}\, H=
(\underline q, \bar q) $, there is a unique root $\beta\in (0,\bar q)$. 

\smallskip
We are looking for conditions under which 
\beq
\label{maina}
0<\liminf_{u\to \infty}u^{\beta}\Psi(u)\le \limsup_{u\to \infty}u^{\beta}\Psi(u)<\infty. 
\eeq
 
 The paper \cite{EKT} treats the case of the non-life insurance.  We formulate its main result in a more  transparent form. 
 \begin{theorem} [\cite{EKT}]
\label{main}
Suppose that the drift $c\ge 0$, the law $F_\xi$  is concentrated on $(-\infty,0)$,   
 $\E[|\xi|^{\beta}]<\infty $,  
 and  $\E[e^{\e T_1}]<\infty$  for some $\e>0$. 
Then (\ref{maina}) holds if  
at least one of the following conditions  are  fulfilled:

 {\bf 1.} $\sigma\neq 0$  or  $\xi$ is unbounded from below. 
 
 {\bf 2a.} $\Pi((-1,0))>0$ and  $\Pi((0,\infty))>0$.
 
 {\bf 2b.}  $\Pi((-1,0))=0$ and  $\Pi(h)=\infty$.  
 
{ \bf 2c.}  $\Pi((0,\infty))=0$ and  $\Pi(|h|)=\infty$. 
 
 {\bf 2d.}    $\Pi((-\infty,0))=0$,
$0<\Pi(h)<\infty$,   
$F((0,t))>0$ for every $t>0$.  

{\bf 2e.}    $\Pi((0,\infty))=0$,
$0<\Pi(|h|)<\infty$,   
$F((0,t))>0$ for every $t>0$.  

\end{theorem}


The proof in \cite{EKT} used heavily  the assumption that  the business process has a  positive drift and  negative claims corresponding to the non-life insurance setting.  In such a case   the ruin may happen only at an instant of jump and, therefore, one needs to monitor the risk process only at $T_1$, $T_2$,   and so on. Such a reduction to a discrete-time ruin model does not work if $\xi_i>0$. 

\smallskip
In our paper we consider the annuity model of the Sparre Andersen type where the ruin occurs because  exhausting resources and the risk process reaches zero in a continuous way. The main result can be formulated as follows. 

 \begin{theorem}
\label{mainAn}
Suppose that the drift $c<0$, the law $F_\xi$   is concentrated on $(0,\infty)$,   
 $\E[\xi^{\beta}]<\infty $,  
 and  $\E[e^{\e T_1}]<\infty$  for some $\e>0$. 
Then (\ref{maina}) holds if  
at least one of the following conditions  are  fulfilled:

 {\bf 1.} $\sigma\neq 0$. 
 
 {\bf 2a.} $\Pi((-1,0))>0$ and  $\Pi((0,\infty))>0$.
 
 {\bf 2b.}  $\Pi((-1,0))=0$ and  $\Pi(h)=\infty$.  
 
{ \bf 2c.}  $\Pi((0,\infty))=0$ and  $\Pi(|h|)=\infty$. 
 
 {\bf 2d.}    $\Pi((-1,0))=0$,
$0<\Pi(h)<\infty$,   
$F((t,\infty))>0$ for every $t>0$.  

{\bf 2e.}    $\Pi((0,\infty))=0$,
$0<\Pi(|h|)<\infty$,   
$F((t,\infty))>0$ for every $t>0$.  

\end{theorem}  

For the mixed case we have the following result.  
\begin{theorem}
\label{mainMix}
Suppose that the drift $c\in \R$,  the law $F_\xi$  charges both  half-lines $(-\infty,0)$ and  $(0,\infty)$,    
 $\E[
|\xi|^{\beta}]<\infty $,  
 and  $\E[e^{\e T_1}]<\infty$  for some $\e>0$. 
Then (\ref{maina}) holds if  
at least one of the following conditions  are  fulfilled:

 {\bf 1.} $\sigma\neq 0$ or  $|\xi|$ is unbounded. 
 
 {\bf 2a.} $\Pi((-1,0))>0$ and  $\Pi((0,\infty))>0$.
 
 {\bf 2b.}  $\Pi((-1,0))=0$ and  $\Pi(h)=\infty$.  
 
{ \bf 2c.}  $\Pi((0,\infty))=0$ and  $\Pi(|h|)=\infty$. 
 
 {\bf 2d.}   $\Pi((-1,0))=0$,
$0<\Pi(h)<\infty$,   
$F((t,\infty))>0$ for every $t>0$ in the case $c<0$ and $F_\xi((0,\e))>0$ for every $\e>0$ in the case $c\ge 0$.

{\bf 2e.}    $\Pi((0,\infty))=0$,
$0<\Pi(|h|)<\infty$,   
$F((t,\infty))>0$ for every $t>0$ in the case $c<0$ and $F_\xi((0,\e))>0$ for every $\e>0$ in the case $c\ge 0$.

\end{theorem}

\smallskip

The annuity and the mixed setting require  a different approach inspired by  the theory of semi-Markov processes. Namely, we consider  
the business process $P$ as the component of the  two-dimensional Markov process $(P,D)$ where the second component  $D=D^r$ is a ``clock", i.e. a process measuring the elapsed time 
after the instant of  last claim.  We  assume that the law of $U=T_1$
may be different from the common law of the further interarrival times:   at the instant zero a portion $r$ of the interarrival time is already elapsed.   This feature admits obvious justifications: e.g., the venture company may change the governance when a project was still in progress.   

Here and throughout the paper we use the superscript $r$  to emphasize  that the law of a random variable or a process depends on $r$, skipping usually $r=0$.   

\smallskip

Formally, the ``clock",  $D^r=(D^r_t)$, is a process with the initial value $D^r_0=r$,   $D^r_t=r+t$ on the interval $[0,T_1)$, and 
  $D^r_t:=t-T^r_n$ on all  other interarrival intervals $[T^r_n,T^r_{n+1})$, 
  $n\ge 1$. That is, the ``clock'' restarts from zero at each instant $T^r_n$. 
  We denote by $F^r$ the law of the first interarrival time 
  $T^r_1=T^r_1-T_0$.  In accordance with our convention $F^0=F$.  
  
  Alternatively, $D^r$ can be representing as the solution of the linear equation 
  $$
  D^r_t=r+t-\int_{[0,t]} D^r_{s-}dN_s. 
  $$

  Typically,
  $\P[T^r_1>t]=\P [T_i>t+r]/\P[T_i>r]$, $i>1$. In the case of exponential distribution $F^r=F$ for all $r\ge 0$ (``absence of memory"). 
 \smallskip

 We  {\bf assume} that $F^r\ge F$.    
 
Recall that the assumed 
independence of $P^r$ and $R$ implies that the joint quadratic characteristic $[P^r,R]$ is zero and  the ruin process $X^{u,r}$ can be written in  the form resembling the Cauchy formula for solutions of liner differential equations:  
\beq
\label{uY}
X^{u,r}_{t}=e^{V_{t}}(u- Y^r_t),
\eeq
where 
\beq
\label{Y_t}
Y^r_t:=-\int_{(0,t]} \cE^{-1}_{s-}(R)dP^r_s=-\int_{(0,t]}  e^{-V_{s-}}dP^r_s.  
\eeq


 The strict positivity of the process $\cE(R)=e^V$ implies that  the ruin time   
 $$
  \tau^{u,r}:=\inf \{t\ge 0:\ X^{u,r}_{t} \le 0\}=\inf \{t\ge 0:\ Y^r_{t} \ge u\}.
$$   


\smallskip 
The crucial element of our study is the following 

\begin{lemma}
\label{G-Paulsen} 
Suppose that  $Y^r_t\to Y^r_{\infty}$ almost surely as $t\to \infty$  where  $Y^r_{\infty}$ is a finite random variable such that   $\bar G(u,r):=\P[Y^r_\infty>u]>0$ for every  $u>0$ and $r\ge 0$. If 
$\bar G_*:=\inf_q \bar G(0,q)>0$, 
then    
\beq
\label{Paulsen}
\bar G(u,r)\le\, \Psi(u,r)=\frac{\bar G(u,r)}{ \E\left[\bar G(X^{u,r}_{\tau^{u,r}},D^r_{ \tau^{u,r}})\, \vert\, \tau^{u,r}<\infty\right]}\le \frac 1{\bar G_*}{\bar G(u,r)}.
\eeq  
\end{lemma}
\noindent
{\sl Proof.}
Let $\tau$ be an arbitrary stopping time with respect to the  filtration $(\cF^{P,D,R}_t)$. 
As we assume that the finite limit $Y^r_\infty$ exists,  the random variable 
$$
Y^r_{\tau,\infty}:=\begin{cases}
-\lim_{N\to \infty } 
\int_{(\tau,\tau+N]}\,e^{-(V_{s-}-V_{\tau})}dP^r_{s},&\tau<\infty, \\
0, & \tau=\infty,
\end{cases} 
$$ 
is well defined.  On the set $\{\tau<\infty\}$
\beq
\label{YX}
Y^r_{\tau,\infty}=e^{V_\tau}(Y^r_{\infty}-Y^r_\tau)=X^{u,r}_{\tau}  +e^{V_\tau}(Y^r_\infty-u). 
\eeq   Let $\zeta$ be a $\cF_{\tau}^{P,D,R}$-measurable random variable.  

Using the strong Markov property, we get that 
\beq
\label{smp}
\P\left[
Y^r_{\tau,\infty}>\zeta, \  \tau<\infty
\right]=\E\left[\bar G(\zeta,D^r_{\tau} )I_{\{ \tau<\infty\}}
\right]
\eeq


Noting that  $\Psi(u,r):=\P\left[\tau^{u,r}<\infty\right ]\ge \P\left[Y^r_{\infty}>u\right]>0$,   we deduce from here using (\ref{YX}) that 
\bean
\bar G(u,r)&=&
\P\left[
Y^r_{\infty}>u,\ \tau^{u,r}<\infty\right]=
\P\left[Y^r_{\tau^{u,r},\infty}>X^{u,r}_{\tau^{u,r}},\ 
\tau^{u,r}<\infty
\right]\\
&= &\P[\tau^{u,r}<\infty] \E\left[\bar G(X^{u,r}_{\tau^{u,r}},D^r_{ \tau^{u,r}})\, \vert\, \tau^{u,r}<\infty\right]\\
&\ge &\P[\tau^{u,r}<\infty] \E\left[\bar G(0,D^r_{ \tau^{u,r}})\, \vert\, \tau^{u,r}<\infty\right] \\
&\ge& \P[\tau^{u,r}<\infty] \inf_q \bar G(0,q)
\eean
and get the result. \fdem

 \smallskip
 
In view of the above lemma the proof of  the main theorem  is reduced to establishing   the existence of finite limits $Y^r_{\infty}$ and finding the asymptotic of the tail of their distributions. 

\smallskip
Let us introduce the notations 
\beq
Q^r_k:=-\int _{(T^r_{k-1},T^r_{k}]}
e^{-(V_{s-}-V_{T^r_{k-1}})}dP^r_s, \qquad M^r_k:=e^{-(V_{T^r_k}-V_{T^r_{k-1}})}.   
\eeq
\begin{lemma}
\label{yr}
The random variables  $Y^r_\infty $ admit the representations 
$$
Y^r_\infty=Q^r_1+M^r_1\tilde Y^r_\infty,
$$ 
where 
\beq
Q^r_1:=-\int _{[0,T^r_1]}
e^{-V_{s-}}dP^r_s, \qquad M^r_1:=e^{-V_{T^r_1}},  
\eeq
$(Q^r_1,M^r_1)$ and $\tilde Y^r_\infty$ are independent, 
and the laws of $\tilde Y^r_\infty$ and  $Y^0_\infty$ coincide. 
\end{lemma}
\noindent
{\sl Proof.}  Note that 
\bean
Y^r_{T^r_n}&=&-\int _{[0,T^r_1]}e^{-V_{s-}}dP^r_s-\sum_{k= 2}^n
e^{V_{T^r_{k-1}}}\int _{(T^r_{k-1},T^r_k]}
e^{-(V_{s-}-V_{T^r_{k-1}})}dP^r_s\\
&=&Q^r_1+M^r_1 \left(Q^r_2+\sum_{k= 3}^n M^r_2...M^r_{k-1}Q^r_k\right),
\eean
where 
the random variable in the parentheses  is independent of 
$(Q^r_1,M^r_1)$
and has the same distribution as $Y^0_{n-1}$. \fdem

\begin{lemma} 
Suppose that  $Y_\infty$ is unbounded from above.  If $c<0$,  then 
$$
\inf_q \bar G(0,q)\ge \E[\bar G(\xi,0)]>0.
$$
If $c\in \R$ and the distribution function $F^r\le F$, then
$$
\inf_q \bar G(0,q)>0.
$$ 
\end{lemma}
\noindent
{\sl Proof.} Using Lemma \ref{yr} we have: 
\bean
\bar G(0,r)&=&\P[Y^r_\infty>0]=\P[Q_1^r/M_1^r+\tilde Y^r_\infty>0]\\
&=&\int \P\left[|c|e^{V_t}\int_{[0,t]}e^{-V_s}ds-\xi_1+ \tilde Y^r_\infty >0\right]F_{T_1^r}(dt)\\
&\ge&  \P[\tilde Y^r_\infty>\xi_1]= \int \P[\tilde Y^r_\infty>x]F_{\xi}(dx)
=\int \bar G(x,0)F_{\xi}(dx)>0
\eean 
since $Y_\infty$  is unbounded. 
\smallskip

The inspection of the proof reveals that the majority of arguments does work with minor changes also for the case where $c$ is of an arbitrary sign and the law $F_\xi$ charges $(-\infty,0)$ and $(0,\infty)$.  In particular, the proof  
that the finite limit $Y_\infty$ exists remains the same. 

Put $f_t:=|\xi_1|+|c| te^{2V^*_{t}}$, where $V^*_t:=\sup_{s\le t}|V_s|$. Then  
$$
|Q_1^r|/M_1^r\le |\xi_1|+|c| e^{V_{T^r_1}}\int_{[0,T_1^r]}e^{-V_s}ds \le f_{{T_1^r}}.
$$
It follows that 
\bean
\bar G(0,r)& =&\P[\tilde Y^r_\infty>-Q_1^r/M_1^r]=\E \bar G(-Q_1^r/M_1^r,0)  \ge \E\bar  G(|Q_1^r|/M_1^r,0) \\
&\ge& \E\int   \bar G(f_t,0)F^r(dt)=-\E\int F^r(t)d\bar G(f_t,0)\\
&\ge& -\E\int F(t)d\bar G(f_t,0)
\ge \E\int  \bar G(f_t,0)F(dt)>0,
\eean
where  we use the  property $F^r\ge F$. Thus, $\inf_r \bar G(0,r)>0$. 
 \fdem

\section{Tails of solutions of distributional equations} 
As a number of results on the  ruin with investments, the proof is based on the implicit renewal theory.  As in \cite{EKT}
we shall use the following formulation combining several  useful facts:    

\begin{theorem} Suppose that  for  some $\beta>0$,
\begin{align}\label{3.3}
\E[M^\beta]=1, \ \ \ \E[M^\beta\,(\ln M)^+]<\infty, \  \  \ \E[|Q|^\beta]<\infty. 
\end{align} 
Let $Y_\infty$ be the solution of the distributional equation $Y_{\infty} \stackrel{d}{=}Q+MY_{\infty}$ and let $\bar G(u):=\P[Y_\infty>u]$. 
Then $\limsup  u^\beta \bar G(u)<\infty$.  If  the random variable  $Y_{\infty}$ is unbounded from above, then  
$\liminf  u^\beta \bar G(u)>0$. 
\end{theorem}

 In the previous section we introduce a process $Y=(Y_t)$.   Assuming that it has at infinity  a  limit $Y_{\infty}$, which is a finite  unbounded from above  random variable, we have proved that it solves the required distributional equation and its tail function gives lower and upper  bounds for the ruin probability.  It remains to check that 
 the hypotheses  of Theorem  \ref{mainAn} ensure the assumed properties and get the result applying the above theorem. We do this in the next sections.  

\section{The existence of the limit  $Y^r_\infty$} 

First, we recall several results from \cite{EKT}. 

\begin{lemma} [ \cite{EKT}, Lemma 2.1]
\label{sup-e}
Let $T> 0$ be a random variable independent of $R$. Suppose that $\E [e^{\e T}]<\infty$ for some $\e>0$. Let $\beta\in (0,\bar q)$ be the root of the equation $H(q)=0$. If $q\in [\beta, \bar q)$ is such that $H(q)\le \e/2$, then   
\beq
\label{sup}
\E\left[ \sup_{s\le T} e^{-qV_{s}}\right]<\infty.
\eeq
\end{lemma} 

\begin{corollary} 
\label{coro}
Suppose that $\E[e^{\e T_1}]<\infty$ for some $\e>0$.
Let 
$$
\widehat Q_1:=\sup_{t\le T_1}|e^{-V_-}\boldsymbol{\cdot} P_t|.
$$ 
If $\E[|\xi_1|^\beta]<\infty$, then $\E [|\widehat Q_1^\beta]<\infty$. 
\end{corollary}

Though the above assertion is a bit more general than  Corollary 2.2 in \cite{EKT}, the proof is exactly the same.  Note also that it does not depend on the sign of $c$ or $\xi_1$ and needs only 
the integrability of $|\xi_1|^\beta$. It implies, in particular, that $\E [|Q_1^\beta|]<\infty$

\begin{lemma} 
\label{limit}
Suppose that $\E[e^{\e T_1}]<\infty$ and $\E[|\xi_1|^{\beta\wedge \e\wedge 1}]<\infty $ for some $\e>0$.  Then $Y_t\to Y_{\infty}$ almost surely as $t\to \infty$  where  $Y_{\infty}$ is a finite random variable. 
\end{lemma}
\noindent
{\sl Proof.}  The convergence a.s. of the sequence $Y_{T_n}$, $n\ge 1$, to a finite r.v. $Y_\infty$ has been proven in Lemma 
4.1 of \cite{EKT} as well as the fact that $\rho:=\E [M_1^p]<1$ for any $p\in (0,\beta\wedge \e\wedge 1)$.  

Put $I_n:=(T_{n-1},T_n]$ and 
$$
\Delta_n:=\sup_{v\in I_n}\left\vert\int_{(T_{n-1},v]}e^{-V_{s-}}dP_s\right\vert=\prod_{i=1}^{n-1}M_i \sup_{v\in I_n}\left\vert\int_{(T_{n-1},v]}e^{-(V_{s-}-V_{T_{n-1}})}dP_s\right\vert.
$$

By virtue of the Borel--Cantelli lemma, to get the announced result  it is sufficient to show that for every $\delta>0$ 
$$
\sum_{n=1}^\infty \P[\Delta_n\ge \delta ]<\infty. 
$$ 
But this is true because the Chebyshev inequality and the Corollary \ref{coro} imply that $\P[\Delta_n\ge \delta ]\le \delta^{-p} \rho^p \E[\widehat Q_1|^p]$. \fdem 

\smallskip
By Lemma \ref{yr} the sequence $Y^r_{T^r_n}$ converges a.s. to $Y^r_\infty$ and the same arguments as above allows us to conclude that $Y^r_{t}$ also converges. 

\section{ When the distribution of $Y_\infty$ is unbounded from above?}  
The question in the title of the section is  studied in \cite{EKT} for the non-life insurance 
case, i.e. when  $c<0$ and $F_\xi((0,\infty))=1$. 
In the present paper we provide  sufficient conditions for the unboundedness from above for all new cases using the techniques developed in the mentioned paper.    
It is based on   the following elementary  observation: if $f:X\times Y\to \R$ is a measurable function, r.v. $\eta$ and $\zeta $ are independent and have the  laws $F_\eta$ and $F_\zeta$, then  the r.v.  $f(\eta,\zeta)$ is unbounded
from above  provided that there exists  a measurable set $X_0\subseteq X$ with  $F_\eta(X_0)>0$ such that the r.v. $f(x,\zeta)$  is unbounded from above for every $x\in X_0$. 
\smallskip

    Let $A_n:=M_1...M_n$, $n\ge 1$,  $A_0:=0$.  
 
 A tractable sufficient condition is give by the following   
\begin{lemma} [\cite{EKT}, Lemma 5.1]
\label{tail}
If there exists $n\ge 1$ such that  the random variables $Q_1$ and $(Q_1+\dots+A_{n-1}Q_n)/A_n$ 
are unbounded from above,  then  $Y_{\infty}$ is  unbounded from above. 
\end{lemma}

It  usually works already with  $n=1$ but sometimes  we need it with $n=2$. A short look at the expressions 
\bea
Q_1&=&-c\int_0^{T_1} e^{-V_r}dr-e^{-V_{T_1}}\xi_1\\
Q_1/A_1&=&-c e^{V_{T_1}}\int_0^{T_1} e^{-V_r}dr-\xi_1,\\
Q_1/A_2 + Q_2/M_2 &= &- ce^{V_{T_2}}\int _0^{T_2}e^{-V_r}dr -\xi_1 e^{V_{T_2} - V_{T_1}} -\xi_2 
\eea
shows that $Y_\infty$ is unbounded from above when $\xi$ is unbounded from below (of course, the latter property is not fulfilled for  the annuity model). 

Using the above sufficient condition of the unboundedness, we examine various cases. 

{\bf 1.} Let $\sigma\neq 0$. In this case the following lemma is helpful: 
\begin{lemma} 
\label{Wiener}
Let $K>0$,   $\sigma\neq 0$ and $0\le s <t$. Then the  random variables 
\beq
\zeta:= Ke^{\sigma W_t}- \int_0^t e^{\sigma W_r}dr , \qquad  \tilde\zeta:= Ke^{\sigma (W_t-W_s)}- e^{\sigma W_t}\int_0^t e^{\sigma W_r}dr, 
\eeq
are  unbounded from below and from above. 
\end{lemma}

The property that $\zeta$ and $\tilde \zeta$ are unbounded from above has been proven in \cite{EKT}, Lemma 5.2. The unboundedness from below can be established by  similar arguments.   It is also clear, that if $K=0$, then $\zeta$ and $\tilde \zeta$ are unbounded from below. 

 \smallskip
 The process  $\bar V:=V-\sigma W$  is independent of the Wiener process  $W$. If $c<0$, then 
 $$
Q_1\ge   |c|
\inf_{r\le T_1} e^{-\bar V_r} \int_0^{T_1}e^{-\sigma W_r}dr - \xi_1e^{-\bar V_{T_1}}e^{-\sigma W_{T_1}}. 
$$
Using the conditioning with respect to $\bar V$, $\xi_1$, $T_1$  and the previous lemma,  we get  that 
$Q_1$ is unbounded from above.  
Since 
$$
Q_1/A_1\ge |c|e^{\bar V_{T_1}}\inf_{r\le T_1}e^{-\bar V_r}e^{\sigma W_{T_1}}\int_0^{T_1} e^{-\sigma W_r}dr-\xi_1,
$$
we conclude in the same way that $Q_1/A_1$ is unbounded from above. If $c\ge 0$, then 
necessary $F_\xi(-\infty,0))>0$ (recall that  we exclude the case $c\ge 0$, $\xi>0$ when the ruin is impossible).  Lemma 5.2 implies that the random variables  $Q_1$ and $Q_1/A_2 + Q_2/M_2$ are unbounded from above. 

{\bf 2.} Let $\sigma=0$ and let $\xi$ be bounded from below. We treat separately several subcases.  
\smallskip

{\bf 2a.} 
$\Pi((-1,0))>0$ and  $\Pi((0,\infty))>0$. Fix $\e>0$ such that    $\Pi((-1, -\e))>0$ and $\Pi((\e,\infty))>0$ and put      
$$
V^{(1)}:=
I_{\{-1<x< -\e\}}\ln (1+x)*\mu +I_{\{x>\e\}}\ln (1+x)*\mu. 
$$ 
Then the processes $V^{(1)}$ and $V^{(2)}:=V-V^{(1)}$ are independent. 

Note that $V^{(1)}$ is the sum of two independent compound Poisson processes with negative and positive jumps, respectively, and the absolute values of jumps are larger than some constant  $c_\e>0$.  
\begin{lemma}
\label{2a}
Let $K>0$, $t>0$. Then the random variable
\beq
\zeta:= Ke^{-V^{(1)}_t} -\int_0^t e^{-V^{(1)}_r}dr
\eeq
is unbounded from above and from below, the random variable 
$$
\widehat \zeta:=e^{-V^{(1)}_t} \int_0^{t} e^{-V^{(1)}_r}dr.
$$
is unbounded from above. 
\end {lemma}
\noindent
{\sl Proof.} The arguments are simple and  we explain only the idea.  One can   consider trajectories where $V^{(1)}$ has a lot of negative jumps 
  in a neighborhood of zero   while all positive jumps are concentrate  in a neighborhood of $t$. 
  Choosing suitable parameters and using the independence of processes with positive and negative jumps we obtain that, with a strictly positive probability, the first term in the definition of $\zeta $ is arbitrary close to zero while the integral is arbitrary large. Thus, $\zeta$ is unbounded from below. Symmetric arguments lead to the conclusion that $\zeta$ is unbounded from above. \fdem

Let $c<0$. Then the following bounds are obvious:
$$
Q_1\ge  |c|
\inf_{r\le T_1} e^{- V^{(2)}_r} \int_0^{T_1}e^{-V^{(1)}_r}dr - \xi_1e^{-\bar V^{(2)}_{T_1}}e^{-V^{(1)}_{T_1}} ,
$$  
$$
Q_1/A_1\ge  |c|e^{V^{(2)}_{T_1}}\inf_{r\le T_1}e^{-V^{(2)}_r}e^{V^{(1)}_{T_1}}\int_0^{T_1} e^{-V^{(1)}_r}dr-\xi_1.
$$
By conditioning with respect to the random variables $V^{(2)}$, $T_1$, $\xi_1$, which are independent of $V^{(1)}$, and using Lemma \ref{2a} we easily obtain that the random variables 
$Q_1$ and $Q_1/A_1$ are unbounded  from above and, by Lemma \ref{tail}, so is $Y_\infty$. 

Let $c\ge 0$. The same arguments as in \cite{EKT} show that $Q_1$ and $Q_1/A_2 + Q_2/M_2$ are unbounded from above on the non-null set $\{\xi_1<0,\;\xi_2<0\}$.

\medskip
{\bf 2b.} $\Pi((-1,0))=0$,  $\Pi(h)=\infty$.

We use the decomposition of $V$  depending  on the choice of $\e\in (0,1)$. Namely, put 
\bea
\label{V1}
V^\e&:=&I_{\{x\le\e\}} h*(\mu-\nu)+I_{\{x\le\e\}}(\ln (1+x)-h)*\mu, \\
\label{V2}
\tilde V^\e&:=&I_{\{x>\e\}} h*(\mu-\nu)+I_{\{x>\e\}}(\ln (1+x)-h)*\mu.
\eea
Note that $V_t=at+V^\e_t+\tilde V^\e_t$ and 
$$
\tilde V^\e=I_{\{x>\e\}} \ln (1+x)*\mu-I_{\{x>\e\}}h*\nu. 
$$

\begin{lemma}
\label{2b}
Let $K>0$, $t>0$. Then the random variables
\beq
\eta:=  \int_0^t e^{-V_r}dr 
- Ke^{-V_t}, \qquad  \eta':= e^{V_t} \int_0^t e^{-V_r}dr 
\eeq
are unbounded from above. 
\end {lemma}
\noindent
{\sl Proof.} Without loss of generality we  assume that $a=0$.  Fix $N>0$ and  choose $\e>0$ small enough to ensure  that $\Pi(I_{\{x>\e\}} h)\ge N$.  
Let  $\Gamma_\e:=\{\sup_{s\le t} |V^\e_s|\le 1\}$. Let denoting by $J^\e$ and $\bar J^\e$ the processes in the lhs of (\ref{V1}).  Using the Doob inequality and the elementary bound $x-\ln (1+x)\le x^2/2$ for $x>0$, we get that  
\bean
\P\left[\sup_{s\le t} |V^\e_s|> 1\right]&\le &\P\left[\sup_{s\le t} |J^\e_s|> 1/2\right]+
  \P\left[|\bar J^\e_t|> 1/2\right] \\
  &\le & 2\E \left[\sup_{s\le t} |J^\e_s|\right]+2\E\left[|\bar J^\e_t|\right]\\
  &\le &
 2 (I_{\{x\le\e\}} h^2*\nu_t)^{1/2} + I_{\{x\le\e\}} h^2*\nu_t\to 0, \quad \e\to 0. 
  \eean
 Thus, the set $\Gamma_\e$ is non-null, at least, for  sufficiently small $\e$, and on this set  
\beq
\label{eta}
\eta\ge \frac 1e \int_0^t e^{-\tilde V^\e_r}dr 
- Ke^{-\tilde V^\e_t+1}.
\eeq
 On the intersection  $\Gamma_\e\cap  
\{I_{\{x>\e\}} h*\mu_{t/2}=0, \ \ln (1+\e)\mu((t/2,t]\times (\e,1] )\ge Nt+1\}$ we have 
$$
\eta\ge \frac 1e \int_0^{t/2} e^{Nr}dr - Ke=\frac 1{eN}(e^{Nt/2}-1)-Ke. 
$$
Due to the independence of $V^\e$ and $\tilde V^\e$ this intersection  is a non-null set. Since $N$ is arbitrary large, the required property of $\eta$ holds. 

The analysis of $\eta'$ follows the same line. At the first stage we replace  $V$ by $\tilde V$  
and compensate the linear decrease of $V$ by a large number of positive jumps on the second half of the interval $[0,t]$.   
 \fdem
\smallskip

Let $c<0$. Then the random variables 
\bean
Q_1&= &|c|\int_0^{T_1} e^{-V_r}dr-e^{-V_{T_1}}\xi_1,\\
Q_1/A_1&=&|c| e^{V_{T_1}}\int_0^{T_1} e^{-V_r}dr-\xi_1
\eean
are unbounded from above by virtue of the lemma. 

Let $c\ge 0$.  The same arguments as in \cite{EKT} lead to the conclusion that  that $Q_1$ and $Q_1/A_2 + Q_2/M_2$ are unbounded from above on the non-null set $\{\xi_1<0,\;\xi_2<0\}$.  

\medskip
{\bf 2c.}   $\Pi((0,\infty))=0$,   
 $\Pi(|h|)=\infty$.  
 
Let $c<0$.  Again we use the decomposition of $V$  depending  on the choice of $\e\in (0,1)$. Put 
\bea
\label{VV1}
V^\e&:=&I_{\{x\ge -\e\}} h*(\mu-\nu)+I_{\{x\ge-\e\}}(\ln (1+x)-h)*\mu, \\
\label{VV2}
\tilde V^\e&:=&I_{\{x<-\e\}} h*(\mu-\nu)+I_{\{x<-\e\}}(\ln (1+x)-h)*\mu,
\eea
 $L:=I_{\{x<-\e\}} \ln (1+x)*\mu$ and $\Pi_\e:=\Pi(I_{\{x<-\e\}}|h|)\uparrow  \infty $ as $\e\to 0$. Then  
$$
\tilde V^\e_t=I_{\{x<-\e\}} \ln (1+x)*\mu_t-I_{\{x<-\e\}}h*\nu=L_t+\Pi_\e t. 
$$
To prove that $Q_1$ is inbounded from above, we argue as follows. As in the previous subcase 
{\bf 2b.} we reduce the problem to checking that the random variable 
$$
\tilde \eta=\int_0^t e^{-\tilde V^\e_r}dr -Ke^{-\tilde V^\e_r}
$$ 
is unbounded from above. 
Let  $t_1:=t_2/2$ where $t_2:=1/\Pi_\e$.  Note that  $t_2\le t$ when $\e>0$ is sufficiently small.  
On the set $\{L_t=L_{t_1}\}$ we have that 
\bean
\tilde \eta  &\ge& (t_2-t_1) e^{|\bar L_{t_1}|}e^{-\Pi_\e t_2}-K
e^{\bar L_{t_1}}e^{-\Pi_\e t} \\
&= &e^{|\bar L_{t_1}|}( 1/(2e\Pi_\e)  -Ke e^{-\Pi_\e t})\ge 
e^{|\bar L_{t_1}|}/(4e\Pi_\e) 
\eean
for sufficiently small $\e$. Since the r.v. $|\bar L_{t_1}|$ is unbounded from above, so is $Q_1$.

On the set $\{ L_{t}=0,\  \xi_1\le K\} $ non-null for any $\e>0$ and sufficiently large $K$, we have the bound 
$$
Q_1/A_1\ge (|c|/\Pi_\e)(e^{\Pi_\e t}-1)-K. 
$$
It follows that $Q_1/A_1$ is unbounded from above. 

The case $c\ge 0$ is treated as in  \cite{EKT}. 

\medskip
{\bf 2d.}    $\Pi((-\infty,0))=0$,
$0<\Pi(h)<\infty$,  
and $F((t,\infty))>0$ for every $t>0$.  

In this subcase $V_t=L_t-bt$,
where $L_t:=\ln (1+x)*\mu_t$ is an increasing process and $b:=\Pi(h)-a$. Note that $b>0$,  
otherwise we get a contradiction with the existence of  $\beta>0$ such that $\ln \E e^{-\beta V_1}=0$. 

Let $c<0$. Take arbitrary  $N>0$. 
Let $s>0$ be the solution of the equation 
$ (|c|/b) (e^{b s}-1)=N$. Take $K$ large enough to ensure that $\P[\xi\le K]>0$. 
Chose $t>s$ such that $F((s,t))>0$.   
On the  non-null set
$$
\{L_s =0, T_1\in (s,t),\  e^{L_t-L_s}\ge Ke^{bt},\  \xi<K\}
$$ 
we have that 
$$
Q_1\ge (|c|/b)(e^{b s}-1) - e^{b t-L_t}\xi\ge N-1. 
$$
Thus, $Q_1$ is unbounded from above. 

To prove that $Q_1/A_1$ is unbounded from above we take $\e>0$  and $K>0$ such that $F((t,\infty))>0$ and  $F_\xi((0,K))>0$. Setting $c_\e:=(|c|/b)(e^{-b \e}-e^{-2b\e})$ we have that on the non-null set 
$$
\{T_1>\e,\ L_{T_1-\e} =0, L_{T_1} \ge \ln ((N+K)/c_\e),  \xi<K\}
$$ 
we have 
$$
Q_1/A_1\ge  (|c|/b)e^{L_{T_1}}e^{-b T_1}(e^{b (T_1-\e)}-1) - \xi\ge c_\e e^{L_t}-K\ge N. 
$$
So, by Lemma \ref{tail} $Y_\infty$ is unbounded. 

If $c\ge 0$, we proceed as in \cite{EKT}, using the assumption that $F$ charges any neighborhood of zero. 
 
\medskip
{\bf 2e.}    $\Pi((0,\infty))=0$,
$0<\Pi(|h|)<\infty$,   
and $F((t,\infty))>0$ for every $t>0$.  

We have again  $V=L_t-bt$ but now the jump process $L$ is decreasing and the constant $b<0$.   

Let $c<0$. 
Fix $N>0$. Let $s,t>0$ be such that $F((s,t))>0$. On the non-null set  
$$
\{T_1\in (s,t),\ |L_{s/2}| \ge N,\ L_{s/2} =L_{t},\  \xi<e^{|L_{s/2}| }\}
$$ 
we have 
$$
Q_1\ge |c|(T_1/2) e^{|L_{(1/2)T_1}| -|b| T_1} - e^{|L_{(1/2)T_1}| -|b| T_1}\xi\ge  |c|(s/2) e^{N-|b| t}-1. 
$$
Since $N$ is arbitrary, $Q_1$ is unbounded from above. 

For any $t>0$ and $K>0$  on the  non-null set  $\{T_1\ge t, \ L_{t}=0,\  \xi\le K\} $ 
$$
Q_1/A_1\le |c/b|(e^{|b|t}-1) -K. 
$$
This implies that $Q_1/A_1$ is unbounded from above. 

If $c\ge 0$, we again proceed as in \cite{EKT}.

\medskip

\section*{Acknowledgements}  
The research is funded by the grant of RSF $n^\circ$ 20-68-47030 ''Econometric and probabilistic methods for the analysis of financial markets with complex structure``. 


\end{document}